\documentclass[11pt]{article}
\usepackage{amssymb,amsmath}

\newcommand{\cl}{\mathop{\rm cl}}
\newcommand{\fr}{\mathop{\rm fr}}
\newcommand{\inter}{\mathop{\rm Int}}

\newcommand{\dist}{\mathop{\rm dist}}

\newcommand{\Id}{\mathop{\rm Id}}

\newcommand{\rank}{\mathop{\rm rank}}

\renewcommand{\Im}{\mathop{\rm Im}}

\numberwithin{equation}{section}
\newtheorem{lemma}{Lemma}[section]
\newtheorem{theorem}{Theorem}[section]

\title{On the unique determination of domains by the condition of
local isometry of boundaries in the relative metrics
\footnote{
(primary); Key words: intrinsic metric,
relative metric of boundary, local isometry of boundaries,
strict convexity}}
\author{Anatoly~P.~Kopylov\footnote{Sobolev Institute of Mathematics,
Acad. Koptyuga pr. 4, and Novosibirsk State University, Pirogova
str., 2, 630090 Novosibirsk, Russia; apkopylov@yahoo.com}}

\begin{document}

\maketitle


\begin{abstract}
The article contains results of recent author's investigations
of rigidity problems of domains in Euclidean spaces
undertaken for a development of a new approach to the classical
problem about the unique determination of bounded closed convex
surfaces~\cite{Po} which is represented by work~\cite{Ko1} in a
sufficiently complete content.

In the article, the full characterization of a plane domain
$U$
with smooth boundary (i.e., the Euclidean boundary
$\fr U$
of
$U$
is a one-dimensional manifold of class
$C^1$
without boundary) that is uniquely determined in the class of
domains in
$\mathbb R^2$
with smooth boundaries by the condition of local isometry of
boundaries in the relative metrics was proved. In the case
where
$U$
is bounded, the necessary and sufficient condition for the
unique determination of the type under consideration in the class
of all bounded plane domains with smooth boundaries is the
convexity of
$U$.
And if
$U$
is unbounded then its unique determination in the class of all
plane domains with smooth boundaries by the condition of local
isometry of boundaries in the relative metrics is equivalent to
its strict convexity.

In the last section, we consider the case of space domains.

The theorem on the unique determination of a strictly convex
domain in
$\mathbb R^n$,
where
$n \ge 2$,
in the class of all
$n$-dimensional
domains by the condition of local isometry of Hausdorff
boundaries in the relative metrics, which is a generalization
of A.~D.~Aleksandrov's theorem on the unique determination of a
strictly convex domain by the condition of (global) isometry of
boundaries in the relative metrics, is proved.

It is also established that in the case of a plane domain
$U$
with nonsmooth boundary and of a three-dimensional domain
$A$
with smooth boundary, the property of domain to be convex is
no longer necessary for their unique determination by the
condition of local isometry of boundaries in the relative
metrics.
\end{abstract}



\section{Introduction}\label{s1}

Let
$\mathcal U$
be a class of domains (i.e., open connected sets) in real
Euclidean
$n$-dimensional
space
$\mathbb R^n$,
where
$n \ge 2$.
We say (see, e.g.,~\cite{Ko1}) that a domain
$U \in \mathcal U$
is uniquely determined in the class
$\mathcal U$
by the relative metric of its (Hausdorff) boundary if each
domain
$V \in \mathcal U$
whose Hausdorff boundary is isometric to the Hausdorff boundary
of the domain
$U$
with respect to the relative metrics is itself isometric to
$U$
(with respect to Euclidean metric).

\textbf{Remark~1.1.} Suppose that
$U$
is a domain in
$\mathbb R^n$
($n \ge 2$)
and
$\rho_U$
is its intrinsic metric. Consider the Hausdorff completion of
the metric space
$(U,\rho_U)$,
i.e., the completion of this space in intrinsic metric
$\rho_U$.
Identifying the points of this completion that correspond to
points of the domain
$U$
with these points themselves and removing them from the
completion, we obtain a metric space
$(\fr_H U,\rho_{\fr_H U,U})$;
the set
$\fr_H U$
of its elements is called the Hausdorff boundary of the domain
$U$,
and
$\rho_{\fr_H U,U}$
is the relative metric on this Hausdorff boundary. The isometry
of the Hausdorff boundaries of domains
$U$
and
$V$
with respect to their relative metrics means the existence of a
surjective isometry
$f: (\fr_H U,\rho_{\fr_H U,U}) \to (\fr_H V,\rho_{\fr_H V,V})$
between these boundaries.

Results of~\cite{Kor1},~\cite{Kor2},~\cite{Kor3} imply, in
particular, that any bounded domain in
$\mathbb R^n$
is uniquely determined by the condition of isometry of boundaries
in the relative metrics. At the same time, according to results
of~\cite{Bor}, a bounded polygonal plane domain
$U$
is uniquely determined by the condition of local isometry of
boundaries in the relative metrics in the class of all such
domains if and only if the domain
$U$
is convex.

\textbf{Remark~1.2.}
Let
$\mathcal M$
be a class of domains in space
$\mathbb R^n$
with
$n \ge 2$.
Following~\cite{Ko1}, we say that a domain
$U \in \mathcal M$
is uniquely determined in the class
$\mathcal M$
by the condition of local isometry of the (Hausdorff) boundaries
of domains in the relative metrics if, for any domain
$V$
belonging to the class
$\mathcal M$,
the local isometry of its Hausdorff boundary to the Hausdorff
boundary of the domain
$U$
with respect to the relative metrics implies the isometry of the
domains
$U$
and
$V$
(with respect to the Euclidean metric). The local isometry in
the relative metrics between the Hausdorff boundaries
$\fr_H U$
and
$\fr_H V$
of the domains
$U$
and
$V$
means the existence of a bijective mapping
$f: \fr_H U \to \fr_H V$
of these boundaries which is a local isometry with respect to
their relative metrics, i.e., a mapping such that, for any
element
$y \in \fr_H U$,
there exists a number
$\varepsilon > 0$
satisfying the following condition: for any two elements
$a$
and
$b$
from the
$\varepsilon$-neighborhood
$Z(y) = \{z \in \fr_H U: \rho_{\fr_H U,U}(z,y) < \varepsilon\}$
of
$y$,
$\rho_{\fr_H U,U}(a,b) = \rho_{\fr_H V,V}(f(a),f(b))$.
It is clear that
$f^{-1}$
is also a local isometry with respect to relative metrics of
boundaries.

In this paper, we continue the study of the unique determination
of domains by the condition of local isometry of their boundaries
in the relative metrics.

It can be divided into two parts.

The first of them (see, Section~2) is mainly devoted to finding
a full description of conditions which is necessary and
sufficient for a plane domain with smooth boundary to be
uniquely determined by the condition of local isometry of their
boundaries in the class of all domains with smooth boundaries
(in the case of a bounded domain, in the class of all bounded
plane domains with smooth boundaries).

The second part is the last section. In this section, we obtain
some new assertions on the unique determination of space
domains with smooth boundaries by the considering in the article
condition. All of these results emphasize the specific character
of our approach to the problems of rigidity of domains in
$\mathbb R^n$.

Note that below
$[a,b] = \{bt + (1-t)a \in \mathbb R^n: 0 \le t \le 1 \}$,
$[a,b[ = \{bt + (1-t)a \in \mathbb R^n: 0 \le t < 1 \}$
($]a,b] = \{bt + (1-t)a \in \mathbb R^n: 0 < t \le 1 \}$)
and
$]a,b[ = \{bt + (1-t)a \in \mathbb R^n: 0 < t <1 \}$
are the segment (closed interval), the half-open interval and
the interval in
$\mathbb R^n$
with endpoints
$a,b \in \mathbb R^n$,
$a \ne b$.
$\inter I$
is the interior of the segment (of the half-open interval)
$I$,
$\inter ]a,b[ = ]a,b[$.
$B(x_0,r) = \{x \in \mathbb R^n: |x - x_0| < r\}$
is the open ball in
$\mathbb R^n$
of radius
$r$
($0 <r < \infty$)
centered at
$x_0 \in \mathbb R^n$.
$\Id_E$
is the identity mapping of a set
$E$:
$\Id_E(x) = x$
for
$x \in E$.

In what follows, paths
$\gamma: [\alpha,\beta] \to \mathbb R^n$,
where
$\alpha,\beta \in \mathbb R$,
are assumed continuous.

\section{The case of plane domains}\label{s2}

The first main result of the article is the following
theorem.

\begin{theorem}\label{t2.1}
Let
$U$
be a domain in
$\mathbb R^2$
with smooth boundary. Then

$(i)$ if
$U$
is bounded{\rm,} then it is uniquely determined in the class of
all bounded plane domains with smooth boundaries by the condition
of local isometry of boundaries in the relative metrics if and
only if this domain is convex;

$(ii)$ if
$U$
is unbounded{\rm,} then the unique determination of
$U$
in the class of all plane domains with smooth boundaries by the
condition of local isometry of boundaries in the relative metrics
is equivalent to the strict convexity of this domain.
\end{theorem}

\textbf{Remark~2.1.}
Let
$U$
be a domain in
$\mathbb R^n$.
As in~\cite{Ko1}, we say that
$U$
has smooth boundary, Lipschitz boundary if the Euclidean
boundary
$\fr U$
of this domain is an
$(n - 1)$-submanifold
of class
$C^1$
(a Lipschitz submanifold) without boundary in
$\mathbb R^n$.
In the case of domain
$U$
with Lipschitz boundary, its Hausdorff boundary
$\fr_H U$
is in natural way identified with Euclidean boundary and metric
$\rho_{\fr U,U}$,
corresponding to Hausdorff metric can be defined in the following
manner:
$$
\rho_{\fr U,U}(x,y) = \liminf_{x' \to x,y' \to y;\, x',y' \in U}
\{\inf[l(\gamma_{x',y',U})]\},
$$
where
$x,y \in \fr U$
and
$\inf[l(\gamma_{x',y',U})]$
is the infimum of lengths
$l(\gamma_{x',y',U})$
of smooth paths
$\gamma_{x',y',U}: [0,1] \to U$
joining
$x'$
and
$y'$
in
$U$.
Recall also that a domain
$U$
is said to be strictly convex if it is convex and the interior
of the segment joining any two points in its closure
$\cl U$
is contained in
$U$.

\begin{lemma}\label{l2.1}
Let
$U$
and
$V$
be two plane domains with smooth boundaries and
$f: \fr U \to \fr V$
be a bijective mapping which is a local isometry of boundaries
of these domains in the relative metrics. Then
$f$
is a {\rm(}global{\rm)} isometry of boundaries
$\fr U$
and
$\fr V$
in their intrinsic metrics.
\end{lemma}

\begin{lemma}\label{l2.2}
Suppose that domains
$U$
and
$V$
and mapping
$f: \fr U \to \fr V$
satisfy to the hypothesis of Lemma~{\rm\ref{l2.1},} moreover{\rm,}
$\fr U$
is bounded. Then the boundary
$\fr V$
of the domain
$V$
is also bounded and the mapping
$f$
has the following property{\rm:} there exists a number
$\varepsilon > 0$
such that
$\rho_{\fr U,U}(a,b) = \rho_{\fr V,V}(f(a),f(b))$
if
$a,b \in \fr U$
and
$\rho_{\fr U,U}(a,b) < \varepsilon$.
\end{lemma}

\begin{lemma}\label{l2.3}
Under hypothesis of Lemma~{\rm\ref{l2.1}} and an additional supposition
that the boundary
$\fr U$
of the domain
$U$
is connected{\rm,} the boundary
$\fr V$
of the domain
$V$
is also connected.
\end{lemma}

The proofs of these lemmas are sufficiently simple. By this
reason, we omit them.
\vskip3mm

\textit{Proof of Theorem~{\rm\ref{t2.1}}.} {\bf Step 1.} Prove the first
part of assertion
$(i)$,
i.e., prove that if
$U$
is a bounded convex plane domain with smooth boundary then it is
uniquely determined in the class of all bounded plane domains with
smooth boundaries by the condition of local isometry of boundaries
in the relative metrics. To this end, suppose that for a bounded
convex plane domain
$U$
with smooth boundary, there exists a bounded plane domain
$V$
with smooth boundary such that its boundary
$\fr V$
is locally isometric to the boundary
$\fr U$
of
$U$
in the relative metrics of boundaries (further,
$f: \fr U \to \fr V$
is the fixed mapping realizing a such isometry). Then by
Lemmas~\ref{l2.2} and~\ref{l2.3}, the boundary of the domain
$V$
is connected (and consequently,
$V$
is a Jordan domain), and
$f$
has the property indicated in Lemma~\ref{l2.2}. Accomplishing, if
it is necessary, an additional inversion with respect to a
straight line, we can also assume that
$f: \fr U \to \fr V$
preserves the orientation of the boundary
$\fr U$
of
$U$,
induced by the canonical orientation of this domain, i.e.,
$f$
"transfers" the mentioned orientation to the orientation of the
boundary
$\fr V$
of
$V$
induced by the canonical orientation of
$V$.

Let, further,
$I = [a,b]$,
where
$a \ne b$,
be a segment such that
$I \subset \fr U$
and the image
$f(I)$
of
$I$
is no longer a segment, moreover, any another segment
$I^* = [a^*,b^*]$
($a^* \ne b^*$,
$I^* \subset \fr U$)
of
$\fr U$
having common points with
$I$
is a subset of the segment
$I$
(below, we denote the set of all such segments
$I$
by the symbol
$\Lambda$).
We assert that
$f(I)$
is a locally convex arc directed by its convexity inside the
domain
$V$.
The latter means that every point
$P \in f(I)$
has a closed neighborhood
$N = N(P)$
for which
$f(I) \cap N$
is a convex arc directed by its convexity inside
$V$,
i.e.,
$f(I) \cap N(P) = f(I_P)$,
where
$I_P = [\alpha_P,\beta_P] \subset I$,
and the closed curve
$C_P$,
composed from
$f(I_P)$
and the segment
$J_P$
joining the endpoints
$f(\alpha_P)$
and
$f(\beta_P)$
of the arc
$f(I_P)$,
either degenerates to the segment
$J_P$,
or is the boundary of a bounded convex domain with the following
property. There is found a segment
$T$
with
$\inter T \subset V$,
placed on straight line
$\tau_P$
which is perpendicular to
$J_P$
and passes through its midpoint, moreover, some endpoint of
$T$
belongs to the arc
$f(I_P)$
and its second endpoint is on the arc
$(\fr V) \setminus f(I_P)$,
both of these endpoints are on the same side from the straight
line
$j_P$
containing the segment
$J_P$,
and the endpoint belonging to the arc
$f(I_P)$
is nearer to
$j_P$
than the other endpoint. Assuming the contrary, i.e., supposing
that
$f(I)$
is no a locally convex arc directed by its convexity inside
$V$,
we (taking the smoothness of the boundary of
$V$
into account)
arrive to the existence of a segment
$I_P = [\alpha_P,\beta_P] \subset I$
such that either
$(1)$
$f(I_P)$
is a nonconvex arc, or
$(2)$
the arc
$f(I_P)$
is no a segment and is a convex arc directed by its convexity
inside the complement
$cV$
of
$V$.
In both cases, for the curve
$f(\inter I_P)$,
there exist a point
$Q \in f(\inter I_P)$
and a locally supporting segment to this curve from the side of
the complement
$cV$
of
$V$
all points of which, except the point
$Q$,
belong to the interior of
$cV$
and
$Q$
is the common point of this segment and the boundary
$\fr V$
of
$V$.
In the case of
$(2)$,
these point and supporting segment can be found on the basis of
the considerations using in the proof of the Theorem of
Leja-Wilkosz~\cite{LW} which is mentioned in~\cite{BZ}, if we
bring evident modifications corresponding to our case in it.

In the case of
$(1)$,
the curve
$f(\inter I_P)$
contains a point
$L_P$
such that if we draw the tangent in
$L_P$
to our curve then there exist points
$R_P \in f(\inter I_P)$
and
$S_P \in f(\inter I_P)$
lying on various sides of this tangent. Replace the point
$R_P$,
if it is necessary, by the point which is the nearest point to
$L_P$
on the segment
$[R_P,L_P]$
(we will remain for this point the designation
$R_P$)
and belongs to the arc
$f(\inter I_P)$.
Analogously, we replace the point
$S_P$
by the point of the arc
$f(\inter I_P)$
which is the nearest point to
$L_P$
on the segment
$[L_P,S_P]$.
And then consider two Jordan domains such that the boundary
of the first of them is the union of the segment
$[R_P,L_P]$
and of that arc from three arcs constituting the set
$f(\inter I_P) \setminus \{R_P,L_P\}$,
the endpoints of which are
$R_P$
and
$L_P$,
and the boundary of the second domain is constructed by the same
way on the basis of the points
$L_P$
and
$S_P$
and of the same arc
$f(\inter I_P)$.
By the way of constructing, one of these domains will be
contained in
$V$
and the second domain will be contained in
$cV$.
Considering the first domain of them and using the above-mentioned
argument from the proof of theorem of Leja-Wilkosz in~\cite{BZ},
it is not difficult to find the desired point
$Q$
on the part of its boundary disposed on
$f(\inter I_P)$,
and a locally supporting segment
$j$
to the curve
$f(\inter I_P)$
at that point from the side of the complement
$cV$
of the domain
$V$.
Hence, in both cases~(1) and~(2), we arrived to desired situation.
Transposing the tangent in the point
$Q$
to
$f(\inter I_P)$
in a parallel way to itself at a sufficiently short distance to
it to the side where, so to say, the domain
$V$
lies, we will easily get the following situation: there exist
three points
$R'_P$,
$L_P$
and
$S'_P$
on the boundary
$\fr V$
of
$V$
belonging to the arc
$f(\inter I_P)$
and such that
$]R'_P,S'_P[ \subset V$,
moreover, the segment
$[R'_P,S'_P]$
cuts off from
$V$
the Jordan subdomain the boundary of which contains
$L_P$.
Clear that
$f^{-1}(R'_P)$,
$f^{-1}(L_P)$
and
$f^{-1}(S'_P)$
are the successively ordered points on the interval
$\inter I_P$.
Hence, the triangle inequality holds for these points in the metric
$\rho_{\fr U,U}$,
but by their choice, for the points
$R'_P$,
$L_P$
and
$S'_P$,
the strict triangle inequality in the metric
$\rho_{\fr V,V}$
holds.
Since we could initially assume that the length of
$I_P$
is less than
$\varepsilon$,
where
$\varepsilon$
is the number from Lemma~\ref{l2.2} corresponding to the mapping
$f$
which is considered now then we arrived to the contradiction
because by this lemma, the equality in the triangle inequality
must also be fulfilled for the points
$R'_P$,
$L_P$
and
$S'_P$.
Therefore,
$f(I)$
is a locally convex arc directed by its convexity inside
$V$.

We assert that the set
$\Lambda$
is finite.
Clear that by the finiteness of the length
$l = l(\fr U)$
of the boundary
$\fr U$
of
$U$,
the finiteness of
$\Lambda$
follows from the fact that
$\Lambda$
does not contain segments the length of which does not exceed,
for example,
$\varepsilon/2$.
Assuming that a segment
$\Delta = [\alpha_{\Delta},\beta_{\Delta}] \in \Lambda$
has the length
$l(\Delta) \le \varepsilon/2$,
consider points
$Q$
and
$S$
of this segment such that
$Q \ne S$,
$Q$
is situated nearer, let us say, to the left endpoint
$\alpha_{\Delta}$
of the segment, and
$f(Q)$
and
$f(S)$
lie on the same side (and at the positive distance) from the
tangent
$\tau$
to
$\fr V$
at the point
$f(\alpha_{\Delta})$,
finally, the (least positive) angle between the tangent rays to
the arcs
$(\fr V) \setminus f(\Delta)$
and
$f([\alpha_{\Delta},S])$
in the points
$f(\alpha_{\Delta})$
and
$f(S)$,
respectively, is less than
$\pi/4$.
Further, let a point
$P \in (\fr U) \setminus \Delta$
is so close to
$\alpha_{\Delta}$
that
$\rho_{\fr U,U}(P,\alpha_{\Delta}) < \varepsilon/2$
and the points
$f(P)$
and
$f(\alpha_{\Delta})$
lie on the same side from each of the tangents to
$\fr V$
in the points
$f(Q)$
and
$f(S)$.
Under these suppositions, for the triple of the points
$P$,
$Q$
and
$S$,
the strict triangle inequality in the metric
$\rho_{\fr U,U}$
holds, and for their images
$f(P)$,
$f(Q)$
and
$f(S)$,
the triangle equality (in the metric
$\rho_{\fr V,V}$)
holds. Thereby, by virtue of the choice of the number
$\varepsilon$
(and Lemma~\ref{l2.2}), we arrive to the contradiction from which
it is follows that
$\Delta = \varnothing$
and consequently,
$\Lambda$
is finite.

Let
$\omega: [0,l] \to \fr U$
be a natural parametrization of the boundary
$\fr U$
of
$U$
corresponding to the orientation of
$\fr U$
generated by the canonical orientation of the domain
$U$,
and let
$[\alpha_1,\beta_1] \subset [0,l]$
and
$[\alpha_2,\beta_2] \subset [0,l]$
be the segments such that
$\omega([\alpha_j,\beta_j]) \in \Lambda$,
where
$j = 1,2$,
$\alpha_1  < \beta_1 < \alpha_2 <\beta_2$
and the arc
$\omega(]\beta_1,\alpha_2[)$
does not contain points of the segments from
$\Lambda$.
We assert that
$f|_{\omega([\beta_1,\alpha_2])}$
is an Euclidean isometry (i.e., there exists an Euclidean isometry
$F: \mathbb R^2 \to \mathbb R^2$
such that
$F|_{\omega([\beta_1,\alpha_2])} = f|_{\omega([\beta_1,\alpha_2])}$).
Indeed, if the arc
$\omega([\beta_1,\alpha_2])$
does not contain segments, i.e., it is strictly convex, moreover,
its convexity directed toward the interior of the complement
$cU$
of
$U$.
Hence, considering a point
$c \in \omega(]\beta_1,\alpha_2[)$
and sufficiently close to it points
$a,b \in \omega(]\beta_1,\alpha_2[)$,
where
$\beta_1 < \omega^{-1}(a) < \omega^{-1}(c) < \omega^{-1}(b) < \alpha_2$
(the closeness of the points
$a$
and
$b$
to the point
$c$
is such that the distance between each pair of the considering
below triple of the points
$a$,
$f^{-1}(\gamma(s_0))$
and
$b$
is less than
$\varepsilon/2$;
we can easily obtain the latter using the hypothesis of theorem)
and assuming that
$[f(a),f(b)] \cap \inter(c V) \ne\varnothing$,
we arrive to a situation where for the shortest path
$\gamma: [0,s] \to \cl V$
joining the points
$f(a)$
and
$f(b)$
in the closure
$\cl V$
of the domain
$V$\footnote{The existence of such shortest path is
guaranteed, for instance, by the results of~\cite{KK1}.},
there exists a point
$s_0 \in ]0,s[$
for which
$\gamma(s_0) \in \fr V$
and
$f^{-1}(\gamma(s_0)) \in f^{-1}(\Im \gamma \cap \fr V) \setminus
\{a,b\}$
($\ne\varnothing$).
But then for the triple of the points
$a$,
$f^{-1}(\gamma(s_0))$
and
$b$,
the strict triangle inequality in the metric
$\rho_{\fr U,U}$
holds,
at the same time for the points
$f(a)$,
$\gamma(s_0)$
and
$f(b)$,
takes place the equality in the triangle inequality in the metric
$\rho_{\fr V,V}$.
Therefore, by virtue of Lemma~\ref{l2.2},
$[f(a),f(b)] \subset \cl V$
from which the equality
$|f(a) - f(b)| = |a - b|$
follows. Hence, the restriction
$f|_{U_{\varepsilon} \cap \omega([\beta_1,\alpha_2])}$
of
$f$
to the intersection
$U_{\varepsilon} \cap \omega([\beta_1,\alpha_2])$
of the
$\varepsilon$-neighborhood
$U_{\varepsilon}$
($= B(P,\varepsilon)$)
of each point
$P \in \omega([\beta_1,\alpha_2])$
and arc
$\omega([\beta_1,\alpha_2])$
itself is an isometry in Euclidean metric. This circumstance
allows easily to conclude that
$f|_{\omega([\beta_1,\alpha_2])}$
is an Euclidean isometry. In the case where
$\omega([\beta_1,\alpha_2])$
contains segments (which no longer belong to the set
$\Lambda$
and, consequently, their images under the mapping
$f$
are also segments), the proof of the fact that
$f|_{\omega([\beta_1,\alpha_2])}$
is an Euclidean isometry is close to the proof of this fact in
the previous case, i.e., in the case of the strict convexity of
$\omega([\beta_1,\alpha_2])$.
The difference in the arguments consists of negligible and
easily reproducible details, and we omit them.

Now, we are able to conclude the proof of the first part of
assertion
$(i)$
of our theorem. If the boundary
$\fr U$
of
$U$
is such that
$\Lambda = \varnothing$,
then the first part of
$(i)$
is proved on the basis of the arguments from the previous item.
In the case of
$\Lambda \ne\varnothing$,
consider a segment
$\Delta \in \Lambda$
and accomplishing appropriate translation and rotation in the
plane
$\mathbb R^2$,
get the situation where the segment
$\Delta$
lies on the ordinate axis, its upper endpoint is the origin, and
the domain
$U$
is situated on the left half-plane. Let
$\gamma: [0,l] \to \fr U$
($\gamma(0) = \gamma(l) = (0,0)$)
be a natural parametrization of the boundary
$\fr U$
of
$U$
corresponding to the orientation of
$\fr U$
generated by the canonical orientation of
$U$.
If
$f|_{\gamma([0,l - l(\Delta)])}$
is an Euclidean isometry then we can assume, without loss of
generality, that
$f|_{\gamma([0,l - l(\Delta)])} = \Id_{\gamma([0,l - l(\Delta)])}$.
Taking yet into account that
$f(\gamma([l - l(\Delta),l])) = f(\Delta)$
is not a segment (because of
$\Delta \in \Lambda$),
we see that
$f(\fr U) = \fr V$
is not a closed curve
i.e.,
$f(\gamma(l)) \ne f(\gamma(0))$.
The obtained contradiction leads us to the conclusion that
$\Lambda = \varnothing$.
Thus, in this case, the first part of
$(i)$
is proved.
Further, assume that
$\Lambda$
consists of
$n$
segments
$[\gamma(\alpha_1),\gamma(\beta_1)]$,
$[\gamma(\alpha_2),\gamma(\beta_2)]$,
$\dots$,
$[\gamma(\alpha_n),\gamma(\beta_n)] = \gamma([l - l(\Delta),l])
= \Delta$,
where
$0 < \alpha_1 < \beta_1 < \alpha_2 < \beta_2 < \dots < \alpha_n <
\beta_n = l$.
Since
$f|_{\gamma([0,\alpha_1])}$
is an Euclidean isometry, we can assume, with loss of generality,
that
$f|_{\gamma([0,\alpha_1])} = \Id_{\gamma([0,\alpha_1])}$.
Then, using the induction argument, it is not difficult to show
that the rotation of the vector
$\omega$,
where
$-\omega$
is the unit tangent vector to the curve
$\gamma([l - l(\Delta),l])$
(i.e., to the segment
$\Delta$) in the point
$\gamma(l)$,
is realized (under the action of the mapping
$f$)
at the angle being equal to the following quantity:
$$
V = -\sum_{k = 1}^n \biggl\{
\sup_{\alpha_k \le t_1 < t_2 < \dots < t_{\varkappa + 1} \le \beta_k}
\sum_{\nu = 1}^{\varkappa} |\theta\gamma(t_{\nu + 1}) -
\theta\gamma(t_{\nu})|\biggr\} \ne 0,
$$
where
$\theta\gamma(t)$
is the unit tangent vector to the curve
$\gamma([t,l])$
in the point
$\gamma(t)$
if
$0 < t < l$
and to the curve
$\gamma([0,l])$
in the point
$\gamma(0) = (0,0)$
when
$t = l$.
If
$|V| < 2\pi$
then
$\omega \ne \mu e_2$,
where
$\mu > 0$
and
$e_2$
is the unit base vector of the ordinate axis.
And if
$|V| \ge 2\pi$
then (since
$f$
preserves the orientation of the boundary) without
self-intersections, the curve
$f(\fr V)$
can not be close. In both cases, we got the contradiction with
the fact that the curve
$\fr V$
is closed and smooth. Hence, the first part of the assertion
$(i)$
of our theorem is completely proved.

{\bf Step 2.} Prove the second part of assertion
$(i)$.
Assuming that
$U$
is a bounded nonconvex plane domain with smooth boundary, we will
show that by the appropriate deformation, we can get another
domain
$V$
whose boundary
$\fr V$
is smooth and locally isometric to the boundary
$\fr U$
of
$U$
in the relative metrics
$\rho_{\fr U,U}$
and
$\rho_{\fr V,V}$
of boundaries, and the domains
$U$
and
$V$
themselves are not isometric each other in Euclidean metric,
i.e., there does not exist an Euclidean isometry
$J: \mathbb R^2 \to \mathbb R^2$
such that
$J(U) = V$.
In the case where the boundary of
$U$
is not connected, a construction of the above-mentioned domain
$V$
realizes by a small permutation of a connected component of the
boundary
$\fr U$
when the location of the other connected components leaves fixed.
And if the boundary of
$U$
is connected, i.e.,
$U$
is a Jordan domain, then we will argue in the following way.
By theorem of Leja-Wilkosz~\cite{LW}, there will be found a
"locally strict supporting outwards" segment
$I$
lying in
$U$
except a single interior point for
$I$,
let us say, point
$P$,
which belongs to
$\fr U$.
Consider a closed disk
$K$
centered at
$P$
and having so small radius
$r$
that the boundary circle of this disk intersects with
$I$
in two points and the interior of one of the half-disks
$K_+$
and
$K_-$
such that
$K_+ \cup K_- = K \setminus I$,
for instance,
$\inter K_-$,
does not contain points of the boundary
$\fr U$
of
$U$.
Let
$u$
and
$v$
are two straight lines which are perpendicular to
$I$,
situated on the various sides of the normal
$n$
to it at the point
$P$,
and sufficiently close to
$n$.
Let us consider the nearest points
$L$
and
$S$
of the sets
$u \cap \fr U$
and
$v \cap \fr U$
to the segment
$I$
and join them by the shortest in
$\cl U$
curve
$\mu$.
Moreover, we regard
$r$
so small that the closure of the arc representing itself lesser
of two arcs, which arise on the boundary
$\fr U$
when we remove the points
$L$
and
$S$
from it, is contained in
$(\inter K_+) \cup \{P\}$
and that (by the smoothness of
$\fr U$)
the curve
$\mu$
is convex, smooth and directed by its convexity toward
$U$.
We can get one of two cases:
$(1)$
$\mu \subset \fr U$,
and
$(2)$
$\mu$
contains segments interior of each of them is a subset of
$U$.
Further, consider (in both cases~(1) and~(2)) the points
$L^*$
and
$S^*$
belonging to
$\lambda \cap \mu$
and chosen by the following way:
$L^*$
and
$L$
lie on the same side of both the straight line containing the
segment
$I$
(moreover, the point
$L^*$
is situated nearer to this straight line than
$L$)
and the straight line
$\psi$
perpendicular to
$I$,
passing through the point
$P$,
besides,
$L^*$
is situated nearer to
$\psi$
than
$L$,
finally, the point
$S^*$
is defined in the similar way in comparison with the location of
$S$.
By symbol
$U^*$,
denote the Jordan domain with the boundary
$((\fr U) \setminus \lambda^*) \cup \mu^*$,
where
$\lambda^*$
and
$\mu^*$
are the subarcs of the arcs
$\lambda$
and
$\mu$
with the endpoints
$L^*$
and
$S^*$.

In the case of
$(1)$,
a necessary deformation of the domain
$U = U^*$
realizes in the sufficiently obvious way and is reduced to a
deformation of the curve
$\mu^*$.
The arc
$\mu^*$
replaces by a convex arc
$\widetilde{\mu^*}$
of the same length, lying in the disk
$K$
and also directed by its convexity in
$U^*$
(more exactly, in the new domain
$\widetilde{U} = \widetilde{U^*}$),
and the arc
$(\fr U^*) \setminus \lambda^* = (\fr U^*) \setminus \mu^*$
leaves fixed, moreover, the closed arc
$\widetilde{\mu^*} \cup \{(\fr U^*) \setminus \mu^*\}$
forms the smooth boundary of the new domain
$\widetilde{U^*}$.
It is not difficult to verify that the boundaries of
$U^*$
and
$\widetilde{U^*}$
are locally isometric in the relative metrics
$\rho_{\fr U^*,U^*}$
and
$\rho_{\fr \widetilde{U^*},\widetilde{U^*}}$
(here, as a local isometry in the relative metrics of the
boundaries of
$U^*$
and
$\widetilde{U^*}$,
we can take the mapping
$f$
of these boundaries which is leaving fixed the arc
$(\fr U^*) \setminus \mu^*$).
Clear also that in the process of the construction of our
deformation, we can get the following situation: it is impossible
to map the domain
$U^*$
onto the domain
$\widetilde{U^*}$
by an Euclidean isometry. Consequently,
$\widetilde{U^*}$
is the desired domain
$V$.

In the case of
$(2)$,
the construction of a new domain
$V$
realizes in the following way. If
$\mu^* = \lambda^*$
then
$V$
is constructed as in the case~(1). In the case where
$\mu^*$
contains segments the interior of which lie in
$U$
and their endpoints belong to
$\fr U$
(denote the set of all such segments by the symbol
$\mathcal M^*$),
we, starting from the domain
$U^*$,
realize first the construction of the domain
$\widetilde{U^*}$
circumscribed in the case~(1), but in addition, we will leave
invariant the length of every segment of
$\mathcal M^*$
under the action of the arising (in the process of the
construction) boundary mapping
$f^*: \fr U^* \to \widetilde{U^*}$.
The latter is possible by the large degree of freedom in the
construction of the curve
$\fr \widetilde {U^*}$
which is submitted by the condition to that the curve
$\mu^*$
satisfies in the case~(2)\footnote{In this connection, see
Lemma~\ref{l4.1}.}. In this case, the final mapping
$f: \fr U \to \fr \widetilde{U}$,
where
$\widetilde{U}$
($= V$)
is a desired new domain, is constructed like this: it leaves
fixed the curve
$(\fr U) \setminus \lambda$
and coincides with
$f^*$
on the set
$N = \mu^* \cap \lambda^*$.
And if the arc
$\chi$
with endpoints
$A$
and
$B$
has not common points with
$(\fr U) \setminus \lambda^*$
and is cut off from
$\fr U$
by a segment from
$\mathcal M^*$,
then we subject this curve to the action of the preserving
orientation Euclidean isometry
$J: \mathbb R^2 \to \mathbb R^2$
such that
$J(A) = f^*(A)$
and
$J(B) = f^*(B)$,
and then put
$f|_{\chi} = J|_{\chi}$.
In this case, the domain
$V$
is the Jordan domain with the boundary
$f(\fr U)$
and by the construction,
$f: \fr U \to \fr V$
is a local isometry of the boundaries of
$U$
and
$V$
in their relative metrics, moreover, the large degree of freedom
in the choice of the above-circumscribed deformation of the
domain
$U$
which still takes place, makes possible to realize this
deformation such that the domains
$U$
and
$V$
are not isometric in the Euclidean metric. So, in both cases
$(1)$
and
$(2)$,
we get the following situation: if
$U$
is nonconvex bounded plane domain with smooth boundary then it is
not uniquely determined in the class of all bounded plane domains
with smooth boundaries by the condition of local isometry of
boundaries in the relative metrics. Consequently, the assertion
$(i)$
of the theorem is completely proved.

{\bf Step 3.} Pass to prove the assertion
$(ii)$.
The fact, that an unbounded strictly convex plane domain
$U$
with smooth boundary is uniquely determined in the class of all
plane domains with smooth boundaries by the condition of local
isometry of boundaries in the relative metrics, can be proved on the
basis of the arguments from the proof of the first part of
assertion
$(i)$.
Considering one more plane domain
$V$
with smooth boundary and assuming that the boundaries of the
domains
$U$
and
$V$
are locally isometric in their relative metrics and modify
negligible the arguments from the proof of the assertion~(i), we
establish that
$\fr U$
and
$\fr V$
are isometric in the Euclidean metric from where the isometry of
the domains
$U$
and
$V$
themselves follows.

{\bf Step 4.} Proving the second part of assertion
$(ii)$,
at first, we make sure that if an unbounded plane domain
$U$
with smooth boundary is not convex then by the same method as we
used in the proof of the second part of assertion
$(i)$,
it can be deformed to a domain
$V$
with smooth boundary, moreover, to such domain that the
boundaries
$\fr U$
and
$\fr V$
are found to be locally isometric in their relative metrics, and
for these domains themselves, there does not exist an Euclidean
isometry
$J: \mathbb R^2 \to \mathbb R^2$
with property
$V = J(U)$.
In the considering case, there exists a deformation of the
boundary
$\fr U$
of the domain
$U$
which does not lead us to the desired result, but the
above-mentioned degree of freedom in a choice of a deformation
makes possible to pass easily over this difficulty.

{\bf Step 5.} Now, let
$U$
be an unbounded plane convex domain with smooth boundary which
is not strictly convex. In this case, a construction of an
above-mentioned domain
$V$
achieves by quite simple methods. Indeed, on the boundary
$\fr U$
of our domain, there exists a segment, let us say,
$I$.
We can assume that any other segment having common points with
$I$
and lying on
$\fr U$,
is a subset of
$I$.
Without loss of generality, we will also suppose that
$I$
is the segment of the abscissa axis with the endpoints
$A = (-2l,0)$
and
$B = (2l,0)$
and the domain
$U$
is found in the lower half-plane. Subject the boundary
$\fr U$
of the domain to the following transformation.
The origin divides the boundary for two curves. The curve from
those curves, which contains the segment with endpoints
$(0,0)$
and
$(0,2l)$,
leaves under this transformation fixed. The segment with
endpoints
$(-l,0)$
and
$(0,0)$
is transformed to the quarter of the circle
$\{(x,y) \in \mathbb R^2: x^2 + (y - \frac{2l}{\pi})^2 =
\frac{4l^2}{\pi^2}\}$
with endpoints
$(0,0)$
and
$P = (-\frac{2l}{\pi},\frac{2l}{\pi})$.
The remaining part of the boundary
$\fr U$
is first subjected to the translation at vector
$((1 - \frac{2}{\pi})l,\frac{2l}{\pi})$
and then to the rotation at angle
$-\frac{\pi}{2}$
with respect to
$P$.
As the final result, we get the curve
$\gamma$
dividing the plane on two unbounded domains. That domain from
them which locally adjoins from below to the segment with
endpoints
$(0,0)$
and
$(0,2l)$,
we will take for a domain
$V$.
Easily to verify that the boundary
$\fr V$
of this domain is locally isometric to the boundary
$\fr U$
of the initial domain
$U$
in the relative metrics
$\rho_{\fr U,U}$
and
$\rho_{\fr V,V}$
of boundaries, and the domains
$U$
and
$V$
themselves are not isometric to each other in Euclidean metric.
Thus, assertion
$(ii)$,
and together with it, Theorem~\ref{t2.1} are completely proved.

In connection with Theorem~\ref{t2.1}, it should be noted that
there exists a bounded plane domain
$U$
with smooth boundary which is not uniquely determined in the
class of all plane domains with smooth boundaries by the
condition of local isometry of boundaries in the relative
metrics~(see~\cite{Sl}).

In the case where the boundary of a domain
$U \subset \mathbb R^2$
is not smooth, Theorem~\ref{t2.1} ceases to be valid.
Really, the following assertion is correct.

\begin{theorem}\label{t2.2}
There exists a bounded plane domain
$U$
with Lipschitz boundary such that it is not convex but{\rm,}
at the same time{\rm,} is uniquely determined in the class of all
plane domains by the condition of local isometry of boundaries in
the relative metrics.
\end{theorem}

\textbf{Remark~2.2.}
Theorem~\ref{t2.2} due to M.~V. Korobkov (see~\cite{Ko2}).
An argument of its proof will be discussed below.

\section{Unique determination of space domains}\label{s3}

Now, consider the case of space domains. But first, remind a
number of notion and facts from~\cite{Ko1} which we use in the
paper below.

\textbf{Definition~3.1.}
The support of an element
$a$
of the Hausdorff boundary
$\fr_H U$
of a domain
$U \subset \mathbb R^n$
($n \ge 2$)
is a point
$a' = a'_a$
of the Euclidean boundary
$\fr U$
to which a Cauchy sequence
$\{x_j\}_{j \in \mathbb N}$
of points
$x_j \in U$
representing the element
$a$
converges in the intrinsic metric
$\rho_U$
of
$U$.

\begin{lemma}\label{l3.1}
Every element
$a \in \fr_H U$
has a unique support
$a'_a$.
\end{lemma}

\begin{lemma}\label{l3.2}
The set of supports
$a'_a$
of elements
$a \in \fr_H U$
is everywhere dense {\rm(}in the Euclidean metric{\rm)} on the
Euclidean boundary
$\fr U$
of
$U$.
\end{lemma}

The mapping
$p_U: \fr_H U \to \fr U$
denotes the transformation of points of the Hausdorff boundary
$\fr_H U$
which puts in the accordance to each element
$a \in \fr_H U$
its support
$a' = a'_a$.
\vskip3mm

Below, we will use the following assertion which is a
generalization of Lemma~3.1 from~\cite{Kor2} to the case of
locally isometric mappings of the boundaries of domains.

\begin{lemma}\label{l3.3}
Let
$U,V$
be domains in space
$\mathbb R^n$
$(n \ge 2)$
such that there exists a bijective mapping
$f: \fr_H U \to \fr_H V$
local isometric in the relative metrics of their Hausdorff
boundaries. Then for every element
$w \in \fr_H U,$
there exists a number
$\varepsilon = \varepsilon_w > 0$
satisfying the following condition{\rm:} for any two elements
$a',b' \in \fr U$
such that
$]a',b'[ \subset U$
and the elements
$a,b \in \fr_H U$
generated by the path
$\gamma(t) = tb' + (1 - t)a',$
$t \in [0,1]$
{\rm(}i.e., generated by the Cauchy sequences in the intrinsic
metric
$\rho_U$
of the domain
$U$
$\{\gamma(1/n)\}_{n = 3,4,\dots}$
and
$\{\gamma(1 - 1/n)\}_{n = 3,4,\dots},$
respectively{\rm),} belong to the
$\varepsilon$-neighborhood
$Z(w) = \{z \in \fr_H U: \rho_{\fr_H U,U}(z,w) < \varepsilon\}$
of the element
$w,$
the relation
$]p_Vf(a),p_Vf(b)[ \subset V$
holds.
\end{lemma}

The proof of this lemma differs from the proof of Lemma~3.1 in
~\cite{Kor2} by negligible modifications, therefore, we omit it.
\vskip3mm

Now, suppose that a considering domain is strictly convex. Then
the following theorem is valid.

\begin{theorem}\label{t3.1}
Let
$n \ge 2$.
If a domain
$U$
in space
$\mathbb R^n$
is strictly convex{\rm,} then it is uniquely determined in the
class of all domains in this space by the condition of local
isometry of boundaries in the relative metrics.
\end{theorem}

\textit{Proof.}
Let
$V$
be a domain such that there exists a bijective mapping
$f: \fr_H U \to \fr_H V$
being a local isometry in the relative metrics of the Hausdorff
boundaries
$\fr_H U$
and
$\fr_H V$
of the domains
$U$
and
$V$.
Assume that
$x$
and
$y$
are points of the Euclidean boundary
$\fr U$
of the domain
$U$
(by the strict convexity of
$U$
and Remark~2.1, we can suppose that
$x$
and
$y$
are simultaneously elements of the Hausdorff boundary
$\fr_H U$
of
$U$).
By Lemma~\ref{l3.3}, each element
$w \in \fr_H U$
has an
$\varepsilon_w$-neighborhood
$Z(w) = \{z \in \fr_H U: \rho_{\fr_H U,U}(z,w) < \varepsilon\}$
with the property: for any points
$a,b \in Z(w)$
the relation
$]p_V f(a),p_V f(b)[ \subset V$
holds (as for
$Z(w)$,
see Lemma~\ref{l3.3}).
From this fact, it follows that the mapping
$\bar{f}: \fr U \to \fr V$
such that
$\bar{f}(x) = p_V f(x)$
if
$x \in \fr U$
is locally isometric in Euclidean metric (i.e., if
$w \in \fr U$,
then for each point
$z \in Z(w)$,
there exist a ball
$B_x = B(x,r_x) \subset \mathbb R^n$
and an isometric mapping
$F_x: \mathbb R^n \to \mathbb R^n$
in the Euclidean metric such that
$F_x|_{(\fr U) \cup B_x} = \bar{f}|_{(\fr U) \cup B_x}$\,).

Let
$\bar{f}(\fr U) = T \subset \fr V$.
We assert that the closure
$\cl T$
of the set
$T$
coincides with the Euclidean boundary
$\fr V$
of
$V$.
Assuming that
$M = ((\fr V) \setminus \cl T) \ne\varnothing$,
consider a point
$z \in M$.
Since
$\cl T$
is a closed set then
$\dist \{z,T\} = \dist \{z, \cl T\} > 0$.
Taking yet into account that by Lemma~\ref{l3.2}, the set of supports
of the Hausdorff boundary of a domain is dense on its Euclidean
boundary, we can assert the existence of an element
$a$
of the Hausdorff boundary
$\fr_H V$
such that its support
$a' = p_V a$
satisfies to the condition
$\dist \{a',T\} = \dist \{a',\cl T\} > 0$.
Let
$\widetilde{a} = f^{-1}(a)$.
We have
$\bar{f}(\widetilde{a}) = p_V f (\widetilde{a}) = p_V(f(f^{-1}(a)
= p_V a = a' \in T$.
Therefore,
$\cl T = \fr V$.

Further, show that the mapping
$\bar{f}$
can be extended to an Euclidean isometry
$F: \mathbb R^n \to \mathbb R^n$
of all space
$\mathbb R^n$.
Indeed, let
$a$
and
$b$
be any two points on the Euclidean boundary
$\fr U$
of
$U$.
We will now establish that
\begin{equation}\label{eq3.1}
|\bar{f}(a) - \bar{f}(b)| = |a - b|.
\end{equation}
To this end, consider a path
$\gamma: [0,1] \to \fr U$
the endpoints of which are
$\gamma(0) = a$
and
$\gamma(1) = b$.
Since
$\bar{f}$
is a local isometry in the Euclidean metric then for each point
$t \in [0,1]$,
we can find a ball
$B_t = B(\bar{f}(\gamma(t)),r_t) \subset \mathbb R^n$
such that there exists an isometric in the Euclidean metric
mapping
$F_t: \mathbb R^n \to \mathbb R^n$
with the property
$F_t|_{(\fr U) \cap B_t} = \bar{f}|_{(\fr U) \cap B_t}$.
By the continuity of the path
$\gamma$,
the sets
$\gamma^{-1}((\fr U) \cap B_t)$,
where
$t \in [0,1]$,
form a covering of the segment
$[0,1]$
which is open with respect to
$[0,1]$.
But then we can extract a finite subcovering
$\{E_s = \gamma^{-1}((\fr U) \cap B_{t_s}), s = 1,\dots,k\}$.
If
$E_{s_1} \cap E_{s_2} \ne\varnothing$
where
$1 \le s_1,s_2 \le k$
then
$(\fr U) \cap B_{t_{s_1}} \cap B_{t_{s_2}} \ne\varnothing$.
Taking into account the strict convexity of the domain
$U$,
we easily conclude that
$F_{t_{s_1}} = F_{t_{s_2}}$.
Thereby, we can assert that there exists the single isometric
in the Euclidean metric mapping
$F: \mathbb R^n \to \mathbb R^n$
such that
$F_s = F$
for all
$s = 1,\dots,k$
and, consequently,
$\bar{f}|_{\Im \gamma} = F|_{\Im \gamma}$.
The latter implies the desired equality~(\ref{eq3.1}). And from
it, by its turn (with regard for the above-stated), the assertion
of the theorem follows.
\vspace{3mm}

\textbf{Remark~3.1.}
Theorem~\ref{t3.1} is a generalization of a theorem of A.~D.~Aleksandrov
about the unique determination of the boundary
$\fr U$
of a strictly convex domain
$U \subset \mathbb R^n$
($n \ge 2$)
by the relative metric
$\rho_{\fr U,U}$
(A.~D.~Aleksandrov's theorem was first published
(with his consent) by V.~A.~Aleksandrov in~\cite{Al}). An
important particular case of this theorem is the following
assertion.

\begin{theorem}\label{t3.2}
Let
$U_1$
be a strictly convex domain in
$\mathbb R^n$.
Assume that
$U_2 \subset \mathbb R^n$
is any domain whose closure is Lipschitz manifold {\rm(}such that
$\fr(\cl U_2) = \fr U_2 \ne\varnothing);$
moreover{\rm,}
$\fr U_1$
and
$\fr U_2$
are isometric {\rm(}globally{\rm)} in their relative metrics
$\rho_{\fr U_1,U_1}$
and
$\rho_{\fr U_2,U_2}$.
Then
$\fr U_1$
and
$\fr U_2$
are isometric in the Euclidean metric of the space
$\mathbb R^n$.
\end{theorem}

\textit{Proof of Theorem~{\rm\ref{t2.2}}.}
Assume that
$U$
is a bounded nonconvex domain in
$\mathbb R^2$
with Lipschitz boundary
$\fr U$,
and there exists such point
$P \in \fr U$
that on the set
$\fr U \setminus \{P\}$,
the domain
$U$
is locally strictly convex, moreover, its convexity directed to
the complement
$c U$
of this domain. We assert that
$U$
is uniquely determined by the condition of local isometry of
boundaries in the relative metrics. The proof of this assertion
realizes by the same scheme and with using the same tools as in
the proof of Theorem~\ref{t3.1}, with certain negligible
modifications. We turn our attention to them briefly.

In the considering now case let
$V$
be a one more domain in
$\mathbb R^2$
whose Hausdorff boundary is locally isometric to the Hausdorff
boundary of
$U$,
let
$f: \fr_H U \to \fr_H V$
be a bijection which is a local isometry in the relative metrics
of the Hausdorff boundaries
$\fr_H U$
and
$\fr_H V$
of
$U$
and
$V$,
finally, let
$T = \bar{f}((\fr U) \setminus \{P\})$
(since the boundary
$\fr U$
of
$U$
is Lipschitz, we, taking into account Remark~2.1, identify
$\fr_H U$
with
$\fr U$).
We assert that in this case, just as in the proof of
Theorem~\ref{t3.1},
$\cl T = \fr V$.
The latter can be proved on the basis of the arguments from the
proof of Theorem~\ref{t3.1}. Nevertheless, by Lemma~\ref{l2.2}
and the infinity of that part of the set of supports of the
Hausdorff boundary
$\fr_H V$
which is contained in the set
$M = (\fr V) \setminus \cl T$,
the indicated there point
$a'$
can be chosen so that
$\alpha = \bar{f}^{-1}(a)\,\, (\ne\varnothing) \subset (\fr U)
\setminus \{P\}$.

The further arguments iterate the arguments used in the proof of
Theorem~\ref{t3.1} almost verbatim. By this reason, we omit them.
\vspace{3mm}

As opposed to that what takes place in the case of domains in
$\mathbb R^2$
(see Theorem~\ref{t2.1}), in the case of space domains, under the
decision of problems on their unique determination by the
condition of local isometry of boundaries in the relative metrics,
the condition of convexity of a considerable domain (as in
Theorem~\ref{t2.2}) ceases to be necessary. In fact, the following
assertion holds.

\begin{theorem}\label{t3.3}
In
$\mathbb R^3,$
there exists a domain
$U$
with smooth boundary such that it is uniquely determined in the
class of all three-dimensional domains with smooth boundaries by
the condition of local isometry of boundaries in the relative
metrics but is not convex.
\end{theorem}

\textit{Proof.}
Let our domain
$U$
be made by the following way.

Consider the arc of cardioid
$$
\theta = \{(x,y,z) \in \mathbb R^3: x^2 + z^2 - \sqrt{x^2 + z^2} +
z = 0,\,\, x^2 + z^2 > 0,\,\, x \ge 0,\,\, y = 0\}.
$$
Leaving it fixed except of the part
$\theta_1$
which is cut out from it by the disk
$\{(x,y,z) \in \mathbb R^3: x^2 + z^2 \le \frac{1}{9},\,\,
y = 0\}$,
replace the arc
$\theta_1$
of the cardioid by the arc of the circle
$\{(x,y,z) \in \mathbb R^3: z = 1 -
\sqrt{\frac{2}{3} - x^2},\,\, 0 \le x \le \frac{\sqrt5}{9},\,\,
y = 0\}$.
It is not difficult to verify that under the rotation of the
curve obtained on this way around the axis
$Oz$
(up to the completed rotation), we obtain the closed smooth
surface being the boundary of a three-dimensional nonconvex
Jordan domain which we will accept for the desired domain
$U$,
establishing further that it is uniquely determined in the class
of all domains in
$\mathbb R^3$
with smooth boundaries by the condition of local isometry of
boundaries in the relative metrics.

So, let
$V \subset \mathbb R^3$
be another domain with smooth boundary and
$f: \fr U \to \fr V$
be a bijective mapping of the boundary
$\fr U$
of
$U$
onto the boundary
$\fr V$
of
$V$
which is a local isometry of the boundaries
$\fr U$
and
$\fr V$
in their relative metrics. Consider the curve
$\theta_0 = \theta \setminus \{(x,y,z) \in \mathbb R^3: x^2 +
z^2 \le 1/4,\,\, y = 0\}$.
Under the rotation around the axis
$Oz$,
this part of cardioid forms a locally strictly convex region
$S$
of the boundary of
$U$
directed by its convexity in the complement
$c U$
of this domain. Applying the same technique as in the proof of
Theorem~\ref{t3.1} and being based on Lemma~\ref{l3.3} in addition,
we first see for ourselves that there exists an isometry
$F: \mathbb R^3 \to \mathbb R^3$
in the Euclidean metric such that
$f|_S = F|_S$.

Without loss of generality, we can assume that
$F = \Id_{\mathbb R^3}$.
Suppose that
$S^*$
is the part of the boundary
$\fr U$
of
$U$
obtained under the rotation of the arc
$\theta^* = \cl(\theta \setminus (\theta_1 \cup \theta_0))$,
and consider the intersection of
$S^*$
with a closed half-plane for which the axis
$Oz$
is the boundary. We can also assume that this intersection is
the curve
$\theta^*$.
Now, we show that any two sufficiently close points
$a$
and
$b$
of this curve (note that the degree of closeness of these points
is determined by Lemma~\ref{l2.2} in application to the mapping
$f$)
cut out from it an arc
$ab$
the image of which under the mapping
$f$
is a plane curve. Indeed, considering the third point
$c$
of the arc
$ab$,
taking into account the local strict convexity of the arc
$\theta^*$
(with respect to the plane domain
$U_{x,z} = U \cap \{(x,y,z) \in \mathbb R^3: y = 0\}$,
moreover, by its convexity in the plane
$\tau_{x,y} = \{(x,y,z) \in \mathbb R^3: y = 0\}$,
this arc is directed to the side of the complement
$\tau_{x,y} \setminus U_{x,z}$
of
$U_{x,z}$),
and applying Lemma~\ref{l3.3} to each pair of the triple of points
$a$,
$c$
and
$b$,
we come to the conclusion that the point
$f(c)$
is on the surface
$\widetilde{S}$
formed by the rotation of the points of the arc
$f(ab) = f(a)f(b)$
around the straight line
$\zeta$
passing through the points
$f(a)$
and
$f(b)$,
and the intersection of
$\widetilde{S}$
with each half-plane, whose boundary is
$\zeta$,
has the same length as the arc
$ab$.
If we suppose that the arc
$f(ab)$
is not plane then its length
will be greater than the length of the arc
$ab$.
The latter contradicts to Lemma~\ref{l2.1}. Hence, the arc
$f(ab)$
is plane. Making arguments close to those which is used in the
proof of the first part of assertion
$(ii)$
of Theorem~\ref{t2.1}, we establish further the existence of an
isometry
$F: \mathbb R^3 \to \mathbb R^3$
in the Euclidean metric such that
$F|_{ab} = f|_{ab}$.
Therefore, the arc
$f(ab)$
(together with the arc
$ab$)
is strictly convex and, consequently, if two planes contain the
arc
$f(ab)$
then they coincide. Extending our last considerations to the arc
$\theta^* \cup \theta_0$,
taking into account the above-said, and using the induction
argument, it is not difficult to establish that the curve
$f(\theta^* \cup \theta_0)$
is contained in the plane
$\tau_{x,y}$,
i.e., in the same plane that the curve
$\theta^* \cup \theta_0$.
Using again the arguments from the proof of Theorem~\ref{t2.1},
we come to the assertion that
$f|_{\theta^* \cup \theta_0} = \Id_{\theta^* \cup \theta_0}$.
Considering the rest intersections of the domain
$U$
with half-planes whose border is axis
$Oz$
and taking into account all above-stated, we obtain as the
result such relation
$$
f|_W = {\Id}_W
$$
where
$W$
is the part of the boundary
$\fr U$
of
$U$
which is obtained by the rotation of the arc
$\theta^* \cup \theta_0$
around the axis
$Oz$.

Assume that
$M = f((\fr U) \setminus W) \cap cV \cap \{(x,y,z) \in
\mathbb R^3: z \ge 2/9\} \ne\varnothing$.
Let
$\alpha > 2/9$
and such that
$$
M_{\alpha} = M \cap \{(x,y,z) \in \mathbb R^3: z = \alpha \}
\ne\varnothing
$$
and
$$
M \cap \{(x,y,z) \in \mathbb R^3: z > \alpha \} = \varnothing.
$$
Suppose that, in
$M_{\alpha}$,
there exists a point
$(\bar{x},\bar{y},\alpha)$
such that
$\bar{x}^2 + \bar{y}^2 > 0$.
Without loss of generality, we can set that
$\bar{x}^2 + \bar{y}^2 = \max\limits_{(x,y,z) \in M_{\alpha}}
(x^2 + y^2)$.
Besides, since
$M_{\alpha} \cap f(W) = \varnothing$
then
$(\bar{x},\bar{y},\alpha) \not\in f(W)$.
Further, let
$\chi = \{\bar{x}(1 + \lambda/\sqrt{\bar{x}^2 + \bar{y}^2})e_1 +
\bar{y}(1 + \lambda/\sqrt{\bar{x}^2 + \bar{y}^2})e_2 +
(\alpha - \lambda t)e_3: \lambda \ge 0 \}$
be a ray outgoing from the point
$P_0 = (\bar{x},\bar{y},\alpha)$,
moreover, the value of the parameter
$t\,\,(>0)$
is so small that this ray intersects
$f((\fr U) \setminus W) \setminus \{P_0\}$
and the distance between
$P_0$
and the nearest point
$P$
of the set
$(f((\fr U) \setminus W) \setminus \{P_0\}) \cap \chi$
to it is lesser than the number
$\varepsilon = \varepsilon_{P_0}$
from Lemma~\ref{l3.3} for the mapping
$f^{-1}$
(in this connection, note that the plane
$\tau_{\alpha} = \{(x,y,z) \in \mathbb R^3: z = \alpha\}$
is supporting to the surface
$f((\fr U) \setminus W)$
and, therefore, is a tangent plane to it in all points
$R \in M_{\alpha}$).
Consequently, by the lemma and the fact that the interval
$]P,P_0[$
is contained in
$V$,
the interval
$]f^{-1}(P),f^{-1}(P_0)[$
must be contained in
$U$.
But this is impossible. Therefore, it remains to consider
the case of
$\bar{x} = \bar{y} = 0$.
And in this case, we also have the contradiction, considering,
for example, the ray
$\{\lambda e_1 + (\alpha - \lambda t)e_3: \lambda \ge 0 \}$
as a desired ray, and further, repeating the arguments used in
the previous case.

We must yet discuss the case
$\alpha = 2/9$.
If
$\dist(M \cap \tau_{2/9},W) >0$
then using the arguments from the previous two cases, we see
that this situation is also impossible. Now, let
$\dist(M \cap \tau_{2/9},W) = 0$.
The stated above facts and the smoothness of the boundaries
$\fr U$
and
$\fr V$
of
$U$
and
$V$
imply the following circumstance: for each point
$z^0 \in M_{2/9}\,\, (= M \cap \tau_{2/9})$,
there exists a number
$\varkappa_0 > 0$
such that any ray emitted from
$z^0$
and intersecting the cone
$K = \{(x,y,z) \in \mathbb R^3: z = \frac{1}{7} +
\frac{5}{63}\sqrt{x^2 + y^2},\,\, \frac{1}{7} \le z \le \frac{2}{9}\}$
in a point contained between the planes
$\tau_{2/9}$
and
$\tau_{2/9 - \varkappa_0}$,
has common points with the surface
$(f((\fr U) \setminus W)) \setminus \{z^0\}$
(here we take into account that the generators of the cone
$K$
pass through the points of the boundary of the manifold
$\cl((\fr U) \setminus W)$,
being tangent in these points to the boundary
$\fr V$
of
$V$).
Choosing as
$z^0$
a point which is so near to
$W$
that the segment
$[z^0,\widetilde z]$
(where
$\widetilde z \in K \cap \tau_{2/9 - \varkappa_0/2}$)
of the ray
$\chi$
emitted from it and intersecting with the circle
$K \cap \tau_{2/9 - \varkappa_0/2}$,
has the least of possible lengths of such segments, consider the
nearest point
$P \in (f((\fr U) \setminus W) \setminus \{z^0\}) \cap \chi$
to the point
$z^0$.
Setting in addition that
$|P - z^0| < \varepsilon_{z^0}$
(where
$\varepsilon_{z^0}$
is a number for the mapping
$f^{-1}$
from Lemma~\ref{l3.3}), we can apply the above-mentioned
arguments to make sure that this case is also impossible. At the
final result, we have the inequality
\begin{equation}\label{eq3.2}
f_3(x,y,z) < \frac{2}{9}
\end{equation}
(where
$f = (f_1,f_2,f_3): \fr U \to \fr V$),
which holds for all points
$(x,y,z) \in (\fr U) \setminus W$.

Consider the bounded open set
$A \in \mathbb R^3$
whose boundary is composed from the sets
$f((\fr U) \setminus W)$
and
$\Xi = \{ (x,y,z) \in \mathbb R^3: x^2 + y^2 \le \frac{5}{81},\,\,
z = \frac{2}{9}\}$.
It is the three-dimensional Jordan domain contained in the
complement to
$V$.
Now, we will prove that the domain
$A$
is convex. Assume by contradiction that this is not valid. Using
the proof of theorem of Leja-Wilkosz~\cite{LW} that is set forth
in~\cite{BZ}, we can assert the existence of three points
$X \in \inter A$,
$Y \in \inter A$
and
$Z \in \inter A$
such that
$[X,Y] \subset \inter A$,
$[Y,Z] \subset \inter A$,
$[X,Z] \not\subset \inter A$,
starting from which and fixing the location of plane
$\tau$
containing these points, we can construct in this plane, for
instance, a locally supporting outwards
$A$
concave arc of ellipse
$\gamma$.
And then changing a location of the point
$Z$
in its small spherical neighborhood, we can obtain a continual
family of locally supporting outwards concave arcs of ellipses.
The plane measure of each part of the boundary
$\fr V$
of
$V$
which is found in one of the indicated plane intersections, can
not be positive, since
$\fr V$
is a smooth bounded surface and, consequently, has a finite area.
Therefore, there exist segments
$[a,b]$
of arbitrary small linear sizes such that
$]a,b[ \subset cA$
and
$a,b \in \fr A$,
moreover, we can also assume that
$a,b \not\in \Xi$.
Hence, we are again found that we have the above-discussed
situation in the process of proving of relation~(\ref{eq3.2})
from which it follows that the domain
$A$
is convex.

As the final result of our arguments for the surfaces
$\cl((\fr U) \setminus W)$
and
$f(\cl((\fr U) \setminus W))$,
we are found themselves in the situation of theorem~2 of
Section 7 from Chapter 3 of monograph of
A.~V.~Pogorelov~\cite{Po}. Using it, we see that these surfaces
are equal. Taking into account the latter and also stated above
facts in the process of proving, we can assert that our theorem
is completely proved.

\section{Appendix}

\begin{lemma}\label{l4.1}
Let
$f_1: [0,a^*] \to \mathbb R$
$(a^* >0)$
be convex downwards strictly increasing smooth function such that
$f_1(0) = f'_1(0) = 0,$
moreover{\rm,} the graph
$\Gamma_1$
of this function contains straight line segments{\rm,} the union
of the set
$\mathcal M$
of all such segments is dense in
$\Gamma_1$
and
$(0,0)$
and
$(a^*,f_1(a^*))$
are limit points for the set of the left endpoints of the
segments from
$\mathcal M$
{\rm(}we assume that the segments
$\Delta \in \mathcal M$
are maximal in such sense that any segment
$\widetilde{\Delta} \subset \Gamma_1$
containing
$\Delta$
coincides with it{\rm)}. Then for each number
$\varepsilon > 0,$
there exists a convex downwards strictly increasing smooth function
$f_2: [0,a^*] \to \mathbb R$
differing from
$f_1$
and having the following properties{\rm:}
$||f_2 - f_1||_{C([0,a^*])} \le \varepsilon,$
$f_2(0) = f'_2(0) = 0,$
$f_2(a^*) = f_1(a^*),$
$f'_2(a^*) = f'_1(a^*)$
and the mapping
$F: \Gamma_1 \to \Gamma_2$
of the graphs of the functions
$f_1$
and
$f_2$
defined by the formula
$$
F: (x,y) \mapsto (\varphi^{-1}(x),f_2(\varphi^{-1}(f^{-1}_1(y)))) \in
\Gamma_2, \quad (x,y) \in \Gamma_1,
$$
where
$\varphi: [0,a^*] \to [0,a^*]$
is the diffeomorphic solution of the functional equation
$$
\int\limits_0^{\varphi(x)}\{1 + [f'_1(\varphi)]^2\}^{1/2} d \varphi =
\int\limits_0^x \{1 + [f'_2(t)]^2\}^{1/2} d t,
\quad 0 \le x \le a^*,
$$
is isometric in the intrinsic metrics of the curves
$\Gamma_1$
and
$\Gamma_2$
which transforms each straight line segment of
$\Gamma_1$
to a straight line segment of
$\Gamma_2$
with the same length.
\end{lemma}

\textit{Proof.}
Let
$x_1$,
$x_2$
and
$x_3$
be three points of the interval
$]0,a^*[$
such that
$x_1 < x_2 < x_3$
and these points are the left endpoints of the segments from the
set
$\mathcal M$
(the choice of the points
$x_1$,
$x_2$,
$x_3$
will be made more precise below). Assume that
$k_1$,
$k_2$,
$k_3$
and
$k_4$
are four real positive numbers. We will choose the function
$f_2$
among functions having the following form:
$$
f_2(x) =
\begin{cases}
k_1 f_1(x), & 0 \le x < x_1; \\
(k_1 - k_2)[f_1(x_1) + f'_1(x_1)(x - x_1)] + k_2 f_1(x),
& x_1 \le x < x_2;\\
\sum\limits_{s = 1}^2 (k_s - k_{s + 1})[f_1(x_s) +
f'_1 (x_s)(x - x_s)] + k_3 f_1(x), & x_2 \le x < x_3; \\
\sum\limits_{s = 1}^3 (k_s - k_{s + 1})[f_1(x_s) + f'_1(x_s)(x - x_s)]
+ k_4 f_1(x), & x_3 \le x \le a^*.
\end{cases}
$$

The equalities
$f_2(a^*) = f_1(a^*)$
and
$f'_2(a^*) = f'_1(a^*)$
leads us to the conditions
\begin{equation}\label{eq4.1}
\sum\limits_{s = 1}^3(k_s - k_{s + 1})[f_1(x_s) +
f'_1(x_s)(a^* - x_s)] + (k_4 - 1)f_1(a^*) = 0
\end{equation}
and
\begin{equation}\label{eq4.2}
\sum\limits_{s = 1}^3(k_s - k_{s + 1}) f'_1(x_s) +
(k_4 - 1) f'_1(a^*) = 0.
\end{equation}

The last condition will be result of the demand
$\varphi(a^*) = a^*$.
And since this demand is equality
$$
\int_0^{a^*}\{1 + [f'_1(t)]^2\}^{1/2} d t = \int_0^{a^*}\{1 +
[f'_2(t)]^2\}^{1/2} d t
$$
then we have
\begin{multline}\label{eq4.3}
\int_0^{a^*}\{1 + [f'_1(t)]^2\}^{1/2} -
\\
\sum\limits_{j = 0}^3 \int\limits_{x_j}^{x_{j + 1}}\biggl\{1 +
\biggl[\sum\limits_{s = 1}^j(k_s - k_{s + 1})f'_1(x_s) +
k_{j + 1}f'_1(t)\biggr]^2\biggr\}^{1/2} d t = 0,
\end{multline}
where
$x_0 = 0$,
$x_4 = a^*$
and
$\sum\limits_{s = 1}^0 \dots = 0$.

The element
$(k_1,k_2,k_3,k_4) = (1,1,1,1) \in \mathbb R^4$
is a solution to the system~(\ref{eq4.1})-(\ref{eq4.3}). At the
same time, by the construction, each straight line segment
$\Delta$
of the curve
$\Gamma_1$
is transformed to a straight line segment of the curve
$\Gamma_2$,
moreover,
$l(F(\Delta)) = l(\Delta)$.
Now, it is sufficient to prove that the rank of the Jacobi matrix
of the left parts of the equalities~(\ref{eq4.1})-(\ref{eq4.3})
calculated with respect to the variables
$k_1$,
$k_2$,
$k_3$
and
$k_4$
in the point
$(1,1,1,1)$
is equal to
$3$
under the successful choice of
$x_1$,
$x_2$
and
$x_3$.

To this end, represent the mentioned matrix in the following
form:
\begin{equation}\label{eq4.4}
N = (A_{j s})_{\begin{array}
{l}
j = 1,2,3 \\
s = 1,2,3,4 \end{array}},
\end{equation}
where
$$
A_{11} = -u_1 - f'_1(x_1) =
-\int_0^{x_1} \frac{[f'_1(t)]^2 d t}{\{1 + [f'_1(t)]^2\}^{1/2}} -
f'_1(x_1) \int_{x_1}^{a^*} \frac{f'_1(t) d t}{\{1 +
[f'_1(t)]^2\}^{1/2}},
$$
$$
A_{12} = - \int_{x_1}^{x_2} \frac{f'_1(t)[f'_1(t) -
f'_1(x_1)] d t}{\{1 + [f'_1(t)]^2\}^{1/2}} -
\int_{x_2}^{a^*} \frac{f'_1(t)[f'_1(x_2) -
f'_1(x_1)] d t}{\{1 + [f'_1(t)]^2\}^{1/2}},
$$
$$
A_{13} = - \int_{x_2}^{x_3} \frac{f'_1(t)[f'_1(t) -
f'_1(x_2)] d t}{\{1 + [f'_1(t)]^2\}^{1/2}} -
\int_{x_3}^{a^*} \frac{f'_1(t)[f'_1(x_3) -
f'_1(x_2)] d t}{\{1 + [f'_1(t)]^2\}^{1/2}},
$$
$$
A_{14} = - \int_{x_3}^{a^*} \frac{f'_1(t)[f'_1(t) -
f'_1(x_3)] d t}{\{1 + [f'_1(t)]^2\}^{1/2}},
$$
$$
A_{21} = f_1(x_1) + (a^* - x_1)f'_1(x_1),
$$
$$
A_{22} = f_1(x_2) - f_1(x_1) + (a^* - x_2)f'_1(x_2) -
(a^* - x_1)f'_1(x_1),
$$
$$
A_{23} = f_1(x_3) - f_1(x_2) + (a^* - x_3)f'_1(x_3) - (a^* -x_2)f'_1(x_2),
$$
$$
A_{24} = f_1(a^*) - f_1(x_3) - (a^* - x_3)f'_1(x_3), \,\,\,
A_{31} = f'_1(x_1),
$$
$$
A_{32} = f'_1(x_2) - f'_1(x_1), \,\,\,
A_{33} = f'_1(x_3) - f'_1(x_2), \,\,\, A_{34} = f'_1(a^*) - f'_1(x_3).
$$

The rank of the matrix~(\ref{eq4.4}) coincides with the rank of
the matrix
$$
\widetilde{N} = \biggl(\sum_{\nu = 1}^s A_{j\nu}\biggr)_{
\begin{array}{l} j = 1,2,3 \\ s = 1,2,3,4 \end{array}},
$$
in which
\begin{multline*}
\sum_{\nu = 1}^2 A_{1\nu} = - u_2 - f'_1(x_2)v_2 =
\\
- \int_0^{x_2} \frac{[f'_1(t)]^2 d t}{\{1 + [f'_1(t)]^2\}^{1/2}}\,\,
-f'_1(x_2) \int_{x_2}^{a^*} \frac{f'_1(t) d t}{\{1 + [f'_1(t)]^2\}^{1/2}},
\end{multline*}
\begin{multline*}
\sum_{\nu = 1}^3 A_{1\nu} = - u_3 - f'_1(x_3)v_3 =
\\
- \int_0^{x_3} \frac{[f'_1(t)]^2 d t}{\{1 + [f'_1(t)]^2\}^{1/2}} \,\, -
f'_1(x_3) \int_{x_3}^{a^*} \frac{f'_1(t) d t}{\{1 +
[f'_1(t)]^2\}^{1/2}},
\end{multline*}
$$
\sum_{\nu =1}^4 A_{1\nu} = - u_4 = - \int_0^{a^*} \frac{[f'_1(t)]^2 d t}
{\{1 + [f'_1(t)]^2\}^{1/2}},
$$
\begin{multline*}
\sum_{\nu = 1}^2 A_{2\nu} = f_1(x_2) + (a^* - x_2)f'_1(x_2),
\,\,\,
\sum_{\nu = 1}^3 A_{2\nu} = f_1(x_3) + (a^* - x_3)f'_1(x_3),
\\
\sum_{\nu = 1}^4 A_{2\nu} = f_1(a^*), \,\,\,
\sum_{\nu = 1}^2 A_{3\nu} = f'_1(x_2), \,\,\, \sum_{\nu = 1}^3
A_{3\nu} = f'_1(x_3), \,\,\, \sum_{\nu = 1}^4 A_{3\nu} = f'_1(a^*).
\end{multline*}

Consider the determinant
\begin{multline*}
\delta_1 = \det \biggl\{\biggl(\sum_{\nu = 1}^s A_{j\nu}\biggr)_{
\begin{array}{l}
j = 2,3 \\ s = 3,4 \end{array}}\biggr\} = [f_1(x_3) +
(a^* - x_3)f'_1(x_3)]f'_1(a^*) -
\\
f_1(a^*)f'_1(x_3) =
f'_1(x_3)f'_1(a^*)\biggl\{\frac{f_1(x_3)}{f'_1(x_3)} + a^* - x_3
- \frac{f_1(a^*)}{f'_1(a^*)}\biggr\}
\end{multline*}
($0 < f'_1(x_3) < f'_1(a^*)$
by the hypothesis of the lemma and the choice of the points
$x_1$,
$x_2$
and
$x_3$).
The second factor in the right part of the last equalities
subject to the following transformations:
\begin{multline}\label{eq4.5}
\frac{\delta_1}{f'_1(a^*) f'_1(x_3)} =
\frac{f_1(x_3)}{f'_1(x_3)} -
\frac{f_1(a^*)}{f'_1(a^*)} +
a^* - x_3 =
\\
- \frac{f_1(a^*) - f_1(x_3)}{f'_1(a^*)} -
f_1(x_3)\biggl(\frac{1}{f'_1(a^*)}
- \frac{1}{f'_1(x_3)}\biggr) +
a^* - x_3
=
\\
- \frac{f'_1(\theta)(a^* - x_3)}{f'_1(a^*)} -
f_1(x_3)\frac{f'_1(x_3) - f'_1(a^*)}{f'_1(a^*)f'_1(x_3)} +
a^* - x_3 =
\\
- \biggl\{\frac{f'_1(\theta) - f'_1(a^*)}{f'_1(x_3) - f'_1(a^*)}
(a^* - x_3) + \frac{f_1(x_3)}{f'_1(x_3)}\biggr\}\frac{f'_1(x_3) -
f'_1(a^*)}{f'_1(a^*)}
\end{multline}
($x_3 < \theta < a^*$).
The convexity of
$f_1$
implies that
$$
\biggl|\frac{f'_1(\theta) - f'_1(a^*)}{f'_1(x_3) -
f'_1(a^*)}\biggr| \le 1.
$$
Therefore, if the condition
\begin{equation}\label{eq4.6}
a^* - x_3 < \frac{f_1(x_3)}{f'_1(x_3)}
\end{equation}
holds then
$\delta \ne 0$
and, consequently,
$\rank \widetilde{N} = \rank N \ge 2$.

Analogously, we can establish that if the point
$x_1$
is fixed and the point
$x_2$
($> x_1$)
is so near to
$x_1$
that the condition
\begin{equation}\label{eq4.7}
x_2 - x_1 < \frac{f_1(x_1)}{f'_1(x_1)}
\end{equation}
holds, then
\begin{multline*}
\delta_2 = \det \biggl\{\biggl(\sum_{\nu = 1}^s A_{j\nu}\biggr)_{
\begin{array}{l} j = 2,3 \\ s = 1,2 \end{array}}\biggr\} = \{f_1(x_1) +
(a^* - x_1)f'_1(x_1)\}f'_1(x_2) -
\\
f'_1(x_1)\{f_1(x_2) +
(a^* - x_2)f'_1(x_2)\} =
f'_1(x_1)f'_1(x_2)\biggl\{\frac{f_1(x_1)}{f'_1(x_1)} -
\frac{f_1(x_2)}{f'_1(x_2)} +
\\
x_2 - x_1\biggr\} \ne 0.
\end{multline*}

If we suppose that the first row of the matrix
$\widetilde{N}$
is a linear combination of two other rows of
$\widetilde{N}$
then we obtain two pairs of the relations
$$
- \frac{u_3}{f'_1(x_3)} - v_3 =
C_1\biggl\{\frac{f_1(x_3)}{f'_1(x_3)} + (a^* - x_3)\biggr\} + C_2,
$$
$$
- \frac{u_4}{f'_1(a^*)} = C_1\frac{f_1(a^*)}{f'_1(a^*)} + C_2
$$
and
$$
- \frac{u_1}{f'_1(x_1)} - v_1 =
C_1\biggl\{\frac{f_1(x_1)}{f'_1(x_1)} + (a^* - x_1)\biggr\} +C_2,
$$
$$
- \frac{u_2}{f'_1(x_2)} - v_2 =
C_1\biggl\{\frac{f_1(x_2)}{f'_1(x_2)} + (a^* - x_2)\biggr\} + C_2,
$$
from which it follows that, under the realization of the
conditions~(\ref{eq4.6}) and~(\ref{eq4.7}),
\begin{multline}\label{eq4.8}
C_1 = \biggl\{- \frac{- u_3}{f'_1(x_3)} + \frac{u_4}{f'_1(a^*)} -
v_3\biggr\}\biggl/\biggl\{\frac{f_1(x_3)}{f'_1(x_3)} -
\frac{f_1(a^*)}{f'_1(a^*)} + (a^* - x_3)\biggr\} =
\\
\biggl\{- \frac{u_2}{f'_1(x_2)} + \frac{u_1}{f'_1(x_1)} - v_2 +
v_1\biggl\}\biggl/\biggl\{\frac{f_1(x_2)}{f'_1(x_2)} -
\frac{f_1(x_1)}{f'_1(x_1)} - (x_2 -x_1)\biggr\},
\end{multline}
moreover,
$C_1$
does not depend on the location of the points
$x_1$,
$x_2$
and
$x_3$.

Further, we have
\begin{multline*}
\frac{u_4}{f'_1(a^*)} - \frac{u_3}{f'_1(x_3)} - v_3 =
\frac{u_4 - u_3}{f'_1(a^*)} + u_3\biggl(\frac{1}{f'_1(a^*)} -
\frac{1}{f'_1(x_3)}\biggr) - v_3 =
\\
\biggl\{- \frac{u_3}{f'_1(a^*)f'_1(x_3)} + O(a^* - x_3)\biggr\}
\{f'_1(a^*) - f'_1(x_3)\}.
\end{multline*}
Note that, here, we used the following estimates:
\begin{multline*}
\biggl|\frac{u_4 - u_3}{f'_1(a^*)} - v_3\biggr| =
\biggl|\frac{1}{f'_1(a^*)}
\int_{x_3}^{a^*} \frac{[f_1(t)]^2 d t}{\{1 +[f'_1(t)]^2\}^{1/2}}\,\,
- \int_{x_3}^{a^*} \frac{f'_1(t) d t}{\{1 +[f'_1(t)]^2\}^{1/2}}\biggr| =
\\
\frac{1}{f'_1(a^*)} \biggl|\int_{x_3}^{a^*}
\frac{f'_1(t)[f'_1(t) - f'_1(a^*)] d t}{\{1 + [f'_1(t)]^2\}^{1/2}}\biggr| \le
\frac{f'_1(a^*) -
f'_1(x_3)}{f'_1(a^*)}
\int_{x_3}^{a^*} \frac{f'_1(t) d t}{\{1 + [f'_1(t)]^2\}^{1/2}}
\\
\le \frac{f'_1(a^*) - f'_1(x_3)}{f'_1(a^*)}
\frac{f'_1(a^*)}{\{1 + [f'_1(a^*)]^2\}^{1/2}}.
\end{multline*}
From these calculations and~(\ref{eq4.5}), we will obtain, as a
result, the equality
\begin{multline}\label{eq4.9}
C_1 = \lim_{x_3 \to a^*}\biggl\{- \frac{u_3}{f'_1(x_3)} +
\frac{u_4}{f'_1(a^*)} - v_3\biggr\} \biggl/
\biggl\{\frac{f_1(x_3)}{f'_1(x_3)} - \frac{f_1(a^*)}{f'_1(a^*)} +
(a^* - x_3)\biggr\} =
\\
- \frac{u_4}{f_1(a^*)} \ne 0.
\end{multline}

By the analogy with~(\ref{eq4.9}) and on the basis
of~(\ref{eq4.8}), we can also to establish that
$$
C_1 = - \frac{u_1}{f_1(x_1)}.
$$
And by the convexity downwards of the function
$f_1$,
\begin{multline*}
\frac{u_1}{f_1(x_1)} = \frac{1}{f_1(x_1)} \int_0^{x_1}
\frac{[f'_1(t)]^2 d t}{\{1 + [f'_1(t)]^2\}^{1/2}} =
\\
\frac{1}{f_1(x_1)} \frac{f'_1(x_1)}{\{1 + [f'_1(x_1)]^2\}^{1/2}}
\int_{\xi}^{x_1} f'_1(t) d t =
\frac{f'_1(x_1)}{\{1 + [f'_1(x_1)]^2\}^{1/2}}
\frac{f_1(x_1) - f_1(\xi)}{f_1(x_1)}
\\
= \frac{f'_1(x_1)}{\{1 + [f'_1(x_1)]^2\}^{1/2}}
\frac{f_1(x_1) - f_1(\xi)}{f_1(x_1) - f_1(0)} \le
\frac{f'_1(x_1)}{\{1 + [f'_1(x_1)]^2\}^{1/2}} \to_{x_1 \to 0} 0
\end{multline*}
($0 < \xi < x_1$),
therefore,
$C_1 = 0$.
As the result of that, we have the contradiction
with~(\ref{eq4.9}). The latter, in turn, leads to the relations
$\rank N = \rank \widetilde{N} = 3$.
The proof of the lemma is completed.

In conclusion, note that the main results of our article were
earlier announced in~\cite{Ko2}.
\vskip3mm

{\bf Acknowledgements}
\vskip3mm

The author was partially supported by the Russian Foundation for
Basic Research (Grant 11-01-00819-a), the Interdisciplinary
Project of the Siberian and Far-Eastern Divisions of the Russian
Academy of Sciences (2012-2014 no. 56), the State Maintenance
Program for the Leading Scientific Schools of the Russian
Federation (Grant NSh-921.2012.1) and the Exchange Program
between the Russian and Polish Academies of Sciences (Project
2014-2016).


\begin{thebibliography}{24}

\bibitem{Ko1}A.~P.~Kopylov,
On the unique determination of domains in Euclidean spaces,
\emph{J. of Math. Sciences}, \textbf{153}, no.6, 869-898 (2008).

\bibitem{Kor1}M.~V.~Korobkov,
Necessary and sufficient conditions for the unique determination
of plane domains,
\emph{Dokl. Math.}, \textbf{76}, 722-723 (2007).

\bibitem{Kor2}M.~V.~Korobkov,
Necessary and sufficient conditions for unique determination of
plane domains,
\emph{Siberian Math. J.}, \textbf{49}, no.3, 436-451 (2008).

\bibitem{Kor3}M.~V.~Korobkov,
A criterion for the unique determination of domains in Euclidean
spaces by the metrics of their boundaries induced by the
intrinsic metrics of the domains, \emph{Siberian Advances in
Mathematics}, \textbf{20}, no.4, 256-284 (2010).

\bibitem{Bor}M.~K.~Borovikova,
On isometry of polygonal domains with boundaries locally
isometric in relative metrics,
\emph{Siberian Math. J.}, \textbf{33}, no.4, 571-580 (1993).

\bibitem{LW}F.~Leja, W.~Wilkosz,
Sur une propri\'{e}t\'{e} des domaines concaves,
\emph{Ann. Soc. Polon. Math.}, \textbf{2}, 222-224 (1924).

\bibitem{BZ}Yu.~D.~Burago, V.~A.~Zalgaller,
Sufficient conditions for convexity, \emph{In: Problems of Global
Geometry}, Leningrad: Nauka (1974). P.~3-53 (Zap. Nauchn. Sem.
LOMI; vol.~45).

\bibitem{KK1}A.~P.~Kopylov and M.~V.~Korobkov,
Rigidity Conditions for the Boundaries of Submanifolds in a
Riemannian Manifold,
\emph{Journal of Siberian Federal University. Mathematics} \&
\emph{Physics},
\textbf{9}(3), 320-331 (2016).

\bibitem{Sl}D.~A.~Slutskiy,
On Two Problems in the Theory of Unique Determination of Domains,
\emph{Vestnik, Quart. J. of Novosibirsk State Univ., Series:
Math., mech. and informatics}, \textbf{11}, no.2, 93-104 (2011).

\bibitem{Ko2}A.~P.~Kopylov,
Unique Determination of Domains by the Condition of Local
Isometry of Boundaries in the Relative Metrics,
\emph{Dokl. Math.}, \textbf{78}, no.2, 746-747 (2008).

\bibitem{Al}V.~A.~Aleksandrov,
Isometry of domains in
$\mathbb R^n$
and relative isometry of their boundaries,
\emph{Siberian Math. J.}, \textbf{25}, no.2, 339-347 (1985).

\bibitem{Po}A.~V.~Pogorelov,
\emph{Extrinsic Geometry of Convex Surfaces}, AMS, Providence (1973).

\end{thebibliography}
\end{document}